\documentclass[12pt]{article}  

\usepackage{ucs} 
\usepackage[utf8x]{inputenc} 
\usepackage[T2A]{fontenc}
\usepackage[russian]{babel}  
 \usepackage{amsfonts,amsmath,amssymb,amscd,amsthm,url}
\usepackage{epsfig}

\title{Универсальный граф для графов ширины разрезания 2}

\author{Н. В. Хорошавкина\footnotemark}
\date{}
\begin{document}
\maketitle
\footnotetext{This paper is prepared under the supervision of Arkadiy Skopenkov and is submitted to the Moscow Mathematical Conference for High-School Students. Readers are invited to send their remarks and reports on this paper to mmks@mccme.ru.

Я благодарю А.Б.Скопенкова за постановку задачи для исследования и помощь в создании данной работы. Я также благодарна рецензенту моей работы за важные замечания.

Supported by the Russian Foundation for Basic Research grant 19-01-00169.}

\begin{abstract}
A graph {\it has cutwidth at most 2} if one can number its vertices by $1,\ldots n$ so that for every $i=1,\ldots,n-1$ there are at most 2 edges $(u,v)$ such that $u\le i<v$.

A characterization of graphs having cutwidth at most 2 in terms of prohibited subgraphs was obtained by Y. Lin, A. Yang.
We present an alternative characterization of such graphs in terms of universal graph.

A {\it chain of $n$ cycles} is a union of $n$ cycles $Z_1,\ldots,Z_n$ having vertices $a_j,b_j\in Z_j$, $j=1,\ldots,n$, such that $Z_i\cap Z_j=\emptyset$ when $|i-j|>1$ and  $Z_{j-1}\cap Z_j=b_{j-1}=a_j$ for every $j=2,\ldots,n$.

{\bf Theorem.} A graph has cutwidth at most 2 if and only if it is a  subgraph of a chain of cycles.
\end{abstract}

В данной работе получена характеризация важного свойства графов. Эта характеризация отлична от известной, ее доказательство проще доказательства известной характеризации.

Граф {\it имеет ширину разрезания не больше 2}, если его вершины можно так занумеровать числами $1,\ldots n$, что для каждого $i=1,\ldots,n-1$ найдется не больше двух ребер $(u,v)$ таких, что $u\le i < v$.

Понятие ширины разрезания активно изучается в теории графов, см. например \cite{LY04} и ссылки в этой статье, \cite{K10}. Графы ширины разрезания не больше 2 -- это в точности те графы, которые допускают непрерывное (или кусочно-линейное) отображение в прямую без трехкратных точек, см. \cite[Problem 3.3.2]{Sk18}, \cite[\S1]{Sk16}.

Характеризация графов ширины разрезания не больше 2 с помощью запрещенных подграфов была получена Ю. Лином, А. Янгом \cite{LY04}.
 Мы охарактеризуем такие графы в терминах универсальных графов. 
Возможно, приводимая ниже характеризация неявно имеется в \cite{LY04}. 
Однако мы не только формулируем ее явно, но и приводим простое прямое доказательство (более простое, чем доказательство другой характеризации в  \cite{LY04}).
 
{\it Цепью циклов} назовем объединение $n$ циклов $Z_1,\ldots,Z_n$, имеющих вершины $a_j,b_j\in Z_j$, $j=1,\ldots,n$, и таких, что $Z_i\cap Z_j=\emptyset$ при $|i-j|>1$, и  $ Z_{j-1}\cap Z_j=b_{j-1}=a_j$ для любого $j=2,\ldots,n$.
\smallskip

{\bf Теорема.} {\it Граф имеет ширину разрезания не больше 2 тогда и только тогда, когда он является подграфом некоторой цепи циклов.}
\smallskip

Из этой теоремы можно получить простое
\smallskip

{\bf Следствие.}
{\it Дерево имеет ширину разрезания не больше 2 тогда и только тогда, когда оно является подграфом некоторого графа, гомеоморфного $G_n$.
Здесь $G_n$ --- граф, имеющий $3n-4$ вершины, $n$ из которых образуют путь, каждая вершина которого, кроме крайних, соединена еще с двумя висячими вершинами.}
\smallskip

Действительно, из теоремы следует, что дерево является подграфом некоторой цепи циклов. При этом в дереве нет ни одного цикла, поэтому каждый цикл $Z_i$ данной цепи можно <<разомкнуть>> на одном двух путей, соединяющих в этом цикле вершины $a_i$ и $b_i$, получив граф, гомеоморфный $G_n$.

\begin{proof}[Доказательство теоремы]
Часть <<тогда>> очевидна, поскольку любая цепь циклов имеет ширину разрезания не больше 2. 

Докажем часть <<только тогда>>.
Пусть граф имеет ширину разрезания не больше 2. 
Рассмотрим любую нумерацию вершин графа, удовлетворяющую определению, числами $\{1,2,\ldots,n\}$. 
Отождествим каждую вершину с ее номером. 

Назовем вершину $i\in\{1,2,\ldots,n\}$ графа  {\it разделяющей}, если никакая пара вершин $u, v$ таких, что $u<i<v$, не соединена ребром, и {\it неразделяющей} в противном случае.
Заметим, что вершины $1$ и $n$ являются разделяющими.

Будем считать, что граф связен, поскольку если каждая из его компонент связности является подграфом некоторой цепи циклов, то это верно и для всего графа.

Докажем по индукции по количеству $k\ge2$ разделяющих вершин, что любой связный граф $G$ ширины разрезания не больше 2, имеющий $k$ разделяющих вершин при некоторой нумерации, удовлетворящей определению ширины разрезания, является подграфом некоторой цепи из $k-1$ цикла, в которой $a_1=1$ и $b_{k-1}=n$.

{\it Докажем базу индукции $k=2$.}
Очевидно, что если в некотором связном графе степени всех вершин не превосходят 2, то этот граф является подграфом цикла. Покажем, что в графе $G$ степени всех вершин не превосходят 2.

Предположим, существует вершина $h$ графа $G$ степени не меньше трех. 
Если $h=1$ или $h=n$, то для $i=1$ или $i=n-1$ соответственно мы получаем противоречие с предположением, что ширина разрезания графа $G$ не превосходит 2, поскольку все ребра, исходящие из вершин 1 или $n$, имеют номера, большие, чем 1, или меньшие, чем $n$, соответственно.

Предположим, вершина $h$ отлична от вершин $1$ и $n$. По принципу Дирихле, среди ребер, исходящих из $h$, найдутся два ребра, $e_1$ и $e_2$, вершины которых, отличные от $h$, одновременно больше (или одновременно меньше), чем $h$. Поскольку в графе $G$ нет других разделяющих вершин, кроме 1 и $n$, то вершина $h$ не является разделяющей, а значит, существует ребро $(a,b)$ такое, что $a < h < b$. Но тогда для трех ребер, $e_1, e_2, (a,b)$, при $i=h$, если ребра $e_1$ и $e_2$ ведут в вершины с большими, чем $h$, номерами, и при $i=h-1$, если ребра $e_1$ и $e_2$ ведут в вершины с меньшими, чем $h$, номерами, мы получаем противоречие с предположением, что ширина разрезания графа $G$ не больше 2.

Значит, степень любой вершины графа $G$ не превосходит 2, и поэтому он является подграфом цикла $Z_1$ с вершинами $a_1=1$, $b_1=n.$

{\it Докажем шаг индукции.} 
Предположим, теорема верна для любого графа, имеющего $k\ge2$ разделяющих вершин при некоторой нумерации вершин этого графа, удовлетворяющей определению ширины разрезания не больше 2. 
Покажем, что тогда она верна для любого графа $G$, имеющего $k+1$ разделяющую вершину при некоторой нумерации вершин этого графа, удовлетворяющей определению ширины разрезания не больше 2.

Обозначим через $x$ разделяющую вершину графа $G$ с минимальным номером, большим 1.
Обозначим через $G_{1,x}$ индуцированный подграф графа $G$ с вершинами $\{1,2,\ldots,x\}$, через $G_{x,n}$ -- индуцированный подграф графа $G$ с вершинами $\{x,x+1,\ldots,n\}$. Покажем, что каждая вершина с номером, меньшим (соответственно, большим), чем $x$, не являющаяся разделяющей в графе $G$, не является разделяющей в графе $G_{1,x}$ (соответственно, в графе $G_{x,n}$). 

Для каждой неразделяющей вершины $h$ графа $G$ существуют вершины $u, v$ такие, что $u<h<v$ и  множество ребер графа $G$ содежит ребро $(u,v)$. Поскольку вершина $x$  --  разделяющая, то $u,v$ одновременно больше или одновременно меньше, чем $x$, а значит, ребро $(u,v)$ содержится ровно в одном из графов $G_{1,x}$ или $G_{x,n}$, и в этом графе вершина $h$ не является разделяющей по определению.

Поскольку $G$ не содержит ребер $(u,v)$ таких, что $u<x<v$, то каждое ребро графа $G$ содержится в графе $G_{1,x}$ или $G_{x,n}$, а значит, объединение графов $G_{1,x}$ и $G_{x,n}$ совпадает с графом $G$. Поскольку графы $G_{1,x}$ и $G_{x,n}$ не содержат разделяющих вершин $n$ и 1 соответственно, то число разделяющих вершин в каждом из них меньше, чем число разделяющих вершин в графе $G$. Отсюда, по предположению индукции, графы $G_{1,x}$ и $G_{x,n}$ являются подграфами некоторых цепей циклов, причем первая из них состоит из одного цикла с вершинами $a_1=1$, $b_1=x$, а вторая -- из $k-2$ циклов с вершинами $a_1=x$, $b_{k-2}=n$. Объединение этих двух цепей циклов также является цепью циклов, состоящей из $k-1$ цикла с вершинами $a_1=1$, $b_{k-1}=n$.

\end{proof}

\end{document}